\documentclass[a4paper,11pt]{article}
\usepackage{latexsym}
\usepackage{amssymb}
\usepackage{enumerate}
\usepackage{amsfonts}
\usepackage{amsmath}
\usepackage[dvipdfmx]{graphicx} 
\usepackage[paper=a4paper,left=30mm,right=20mm,top=25mm,bottom=30mm]{geometry}

\newenvironment{@abssec}[1]{%
    \if@twocolumn

      \section*{#1}%
    \else

      \vspace{.05in}\footnotesize
      \parindent .2in
 {\upshape\bfseries #1. }\ignorespaces
    \fi}

    {\if@twocolumn\else\par\vspace{.1in}\fi}

\newenvironment{keywords}{\begin{@abssec}{\keywordsname}}{\end{@abssec}}

\newenvironment{AMS}{\begin{@abssec}{\AMSname}}{\end{@abssec}}

\newcommand\keywordsname{Key words}
\newcommand\AMSname{AMS subject classifications}
\newcommand\AMname{AMS subject classification}
\newtheorem{theorem}{Theorem}

\title{A symmetry theorem in two-phase heat conductors\thanks{This research was partially supported by the Grants-in-Aid
for Scientific Research (B) and (C) ($\sharp$ 18H01126 and $\sharp$ 22K03381)  of
Japan Society for the Promotion of Science and National Research Foundation of S. Korea grant 2022R1A2B5B01001445.}}

\author{Hyeonbae Kang\thanks{Department of Mathematics and Institute of Applied Mathematics, Inha University, Incheon 22212, S. Korea (hbkang@inha.ac.kr).}  \and Shigeru Sakaguchi\thanks{Graduate School of Information Sciences, Tohoku University, Sendai, 980-8579, Japan (sigersak@tohoku.ac.jp).}}

\date{}
\begin{document}
\maketitle

\begin{abstract}
We consider the Cauchy problem for the heat diffusion equation in the whole Euclidean space consisting of two media with different constant conductivities, where initially one medium  has temperature 0 and the other  has temperature 1. Under the assumptions that one medium is bounded and the interface is of class $C^{2,\alpha}$,
we show that if the interface is  stationary isothermic, then it must be a sphere. The method of moving planes due to Serrin is directly utilized to prove the result.
  \end{abstract}



\pagestyle{plain}
\thispagestyle{plain}


\section{Introduction}
\label{introduction}

In the previous paper \cite{KS2021}, we considered the Cauchy problem for the heat diffusion equation in the whole Euclidean space consisting of two media with different constant conductivities, where initially one medium has temperature 0 and the other  has temperature 1. There, the large time behavior, either stabilization to a constant or oscillation,  of temperature was studied.
The present paper deals with the case where  one medium is bounded and the interface is of class $C^{2,\alpha}$,
and  introduces an overdetermined problem with the condition that the interface is  stationary isothermic.

To be precise, let $\Omega$ consist of a finite number, say $m$, of  bounded domains $\{ \Omega_j\}$  in $\mathbb R^N$ with $N \ge 2$, where each $\partial\Omega_j$ is of class $C^{2,\alpha}$ for some $0< \alpha <1$ and  $\overline{\Omega}_i\cap \overline{\Omega}_j=\emptyset$ if $i\not=j$.
Denote by $\sigma=\sigma(x)\ (x \in \mathbb R^N)$  the conductivity distribution of the whole medium given by
\begin{equation}
\label{conductivity constants}
\sigma =
\begin{cases}
\sigma_+ \quad&\mbox{in } \Omega=\bigcup\limits_{j=1}^m\Omega_j, \\
\sigma_- \quad &\mbox{in } \mathbb R^N \setminus \Omega,
\end{cases}
\end{equation}
where $\sigma_-, \sigma_+$ are positive constants with $\sigma_-  \not= \sigma_+$.
The diffusion over such multiphase heat conductors has been dealt with also in \cite{Satrieste2016, SaBessatsu2020, SaJMPA2020, CSU2019, CMS2021}.

We consider the unique bounded solution  $u =u(x,t)$ of the Cauchy problem for the heat diffusion equation:
\begin{equation}
  u_t =\mbox{ div}(\sigma \nabla u)\quad\mbox{ in }\  \mathbb R^N\times (0,+\infty) \ \mbox{ and }\ u\ ={\mathcal X}_{\Omega}\ \mbox{ on } \mathbb R^N\times
\{0\},\label{heat Cauchy}
\end{equation}
where ${\mathcal X}_{\Omega}$ denotes the characteristic function of the set $\Omega$.
The maximum  principle gives
\begin{equation}
\label{positive values}
0 < u(x, t) < 1\ \mbox{ for every } (x,t) \in \mathbb R^N\times (0, + \infty).
\end{equation}

Our symmetry theorem is stated as follows.

\begin{theorem}
\label{symmetry theorem}
If there exists a function $a:(0,+\infty)\to(0,+\infty)$ satisfying
\begin{equation}
\label{stationary isothermic}
u(x,t)=a(t)\ \mbox{ for every } (x,t) \in \partial\Omega\times(0,+\infty),
\end{equation}
then $\Omega$ must be a ball.
\end{theorem}

If $\partial\Omega$ is of class $C^6$, then Theorem \ref{symmetry theorem}
can be proved by the method employed in \cite[Theorem 1.5 with the proof, pp. 335--341]{CMS2021}, where concentric balls are characterized.  The proof there consists of four steps summarized as follows:  (i) reduction of \eqref{heat Cauchy} to elliptic problems by  the Laplace-Stieltjes transform $\lambda\int_0^\infty e^{-\lambda t}u(x,t)dt$ for all sufficiently large $\lambda > 0$,  (ii) construction of precise barriers based on the formal WKB approximation where the fourth derivatives of the distance function to $\partial\Omega$ together with the assumption \eqref{stationary isothermic} are used, (iii) showing that the mean curvature of $\partial\Omega$  is constant with the aid of the precise asymptotics as $\lambda \to \infty$ and the transmission conditions on the interface $\partial\Omega$, (iv) Alexandrov's soap bubble theorem \cite{Al1958} from which we conclude that $\partial\Omega$ must be a sphere.

The approach of the present paper is different from that in \cite{CMS2021} and only requires $\partial\Omega$ to be of class $C^{2,\alpha}$ for some $\alpha>0$. Here the proof consists of  two ingredients:  (i) reduction to elliptic problems by  the Laplace-Stieltjes transform $\lambda\int_0^\infty e^{-\lambda t}u(x,t)dt$ for some $\lambda$, for instance $\lambda =1$,  (ii) the method of moving planes due to Serrin \cite{Se1971, GNN1979, R1997, SiPoincare2001} with the aid of the transmission conditions on $\partial\Omega$. To apply the method of moving planes, the solutions need to be of class $C^2$ up to the interface $\partial\Omega$ from each side, which is guaranteed if $\partial\Omega$ is of class $C^{2,\alpha}$.

\setcounter{equation}{0}
\setcounter{theorem}{0}

\section{Introducing a Laplace-Stieltjes transform}
\label{section_Laplace}

Let $u=u(x,t)$ be the unique bounded solution  of \eqref{heat Cauchy} satisfying \eqref{stationary isothermic}.  We use the Gaussian bounds for the fundamental solutions of diffusion equations due to
Aronson \cite[Theorem 1, p. 891]{Ar1967}(see also \cite[p. 328]{FS1986}). Let $g = g(x,\xi,t)$ be the fundamental solution of $u_t=\mbox{ div}(\sigma\nabla u)$. Then there exist two positive constants $\lambda < \Lambda$ such that
\begin{equation}
\label{Gaussian bounds}
\lambda t^{-\frac N2}e^{-\frac{|x-\xi|^2}{\lambda t}}\le g(x,\xi,t) \le \Lambda t^{-\frac N2}e^{-\frac{|x-\xi|^2}{\Lambda t}}
\end{equation}
 for all $(x,t), (\xi,t) \in \mathbb R^N\times(0,+\infty)$.  Note that $u$ is represented as
 \begin{equation}
 \label{expression of the solution}
 u(x,t)=\int\limits_\Omega g(x,\xi,t)d\xi\ \mbox{ for } (x,t) \in \mathbb R^N\times(0, +\infty).
 \end{equation}

Define the function $v=v(x)$ by
\begin{equation}
\label{Laplace-Stieltjes transform for 1}
v(x) = \int_0^\infty e^{-t}u(x,t) dt\ \mbox{ for } x \in \mathbb R^N.
\end{equation}
With the function $a$ in \eqref{stationary isothermic}, we set $a^*=\int_0^\infty e^{-t}a(t) dt$. Then, \eqref{positive values} yields that $0<a^* < 1$.
Set
 \begin{equation}
 \label{inside and outside}
 v^+=v\ \mbox{ for }x \in \overline{\Omega} \quad \mbox{and} \quad  v^-=v\ \mbox{ for }x \in \mathbb R^N\setminus\Omega.
 \end{equation}
 Then we observe that
\begin{eqnarray}
&&a^* < v^+ < 1 \quad \mbox{and} \quad -\sigma_+\Delta v^++v^+=1\ \mbox{ in } \Omega, \label{interior}\\
&& 0< v^- < a^*  \quad \mbox{and} \quad  -\sigma_-\Delta v^-+v^-=0 \ \mbox{ in } \mathbb R^N\setminus\overline{\Omega},\label{exterior}\\
&& v^+=v^- = a^* \quad \mbox{and} \quad  \sigma_+\frac {\partial v^+}{\partial \nu\ }=  \sigma_-\frac {\partial v^-}{\partial \nu\ } \ \mbox{ on } \partial\Omega,\label{transmission}\\
&& \lim\limits_{|x|\to\infty} v^-(x) = 0.\label{at infinity}
\end{eqnarray}
Here,  $\nu$ denotes the  outward unit normal vector to $\partial\Omega$, the inequalities in \eqref{interior} and \eqref{exterior} follow from the maximum principle, \eqref{transmission} expresses the transmission conditions on the interface $\partial\Omega$, and \eqref{at infinity} follows from \eqref{Gaussian bounds} and \eqref{expression of the solution}.

\setcounter{equation}{0}
\setcounter{theorem}{0}

\section{Proof of Theorem \ref{symmetry theorem}}
\label{section_Proof}

Let us apply  directly the method of moving planes due to Serrin \cite{Se1971, GNN1979, R1997, SiPoincare2001} to our problem in order to show that $\Omega$ must be a ball.
The point is to apply the method to both the interior $\Omega$ and the exterior $\mathbb R^N\setminus\overline{\Omega}$ at the same time. For the method of moving planes for $\mathbb R^N\setminus\overline{\Omega}$, we refer to \cite{R1997, SiPoincare2001}. In this procedure, the supposition that $\Omega$ is not symmetric will lead us to the contradiction that the transmission conditions \eqref{transmission} do not hold.

Let $\gamma$ be a unit vector in $\mathbb R^N,$ $\lambda\in\mathbb R,$ and let $\pi_\lambda$ be the hyperplane $x\cdot\gamma=\lambda.$
For large $\lambda,$ $\pi_\lambda$ is disjoint from $\overline{\Omega}$; as $\lambda$ decreases,
$\pi_\lambda$ intersects $\overline{\Omega}$ and cuts off from $\Omega$ an open cap $\Omega_\lambda= \Omega \cap \{ x \in \mathbb R^N : x\cdot\gamma > \lambda\}$.

Denote by $\Omega^\lambda$ the reflection
of $\Omega_\lambda$ with respect to the plane $\pi_\lambda$. Then, $\Omega^\lambda$ is contained in $\Omega $ at the beginning, and remains in $\Omega$ until one of the following events occurs:
\begin{enumerate}
\item[(i)] $\Omega^\lambda$ becomes internally tangent to $\partial \Omega$
at some point $p\in\partial\Omega\setminus \pi_\lambda;$
\item[(ii)] $\pi_\lambda$ reaches a position where it is orthogonal to $\partial \Omega$
at some point $q\in\partial\Omega\cap\pi_\lambda$ and the direction $\gamma$ is not tangential to $\partial\Omega$ at every point on $\partial \Omega \cap \{ x \in \mathbb R^N : x\cdot\gamma > \lambda\}$.
\end{enumerate}
Let $\lambda_*$ denote the value of $\lambda$ at which either (i) or (ii) occurs.
We claim that $\Omega$ is symmetric with respect to $\pi_{\lambda_*}$.

Suppose that $\Omega$ is not symmetric with respect to $\pi_{\lambda_*}$.
Denote by  $D$  the reflection
of $\left(\mathbb R^N\setminus\overline{\Omega}\right)\cap\{ x \in \mathbb R^N : x\cdot\gamma > \lambda_*\}$ with respect to $\pi_{\lambda_*}$.
Let $\Sigma$ be the connected component of $\left(\mathbb R^N\setminus\overline{\Omega}\right)\cap\{ x \in \mathbb R^N : x\cdot\gamma < \lambda_*\}$ whose boundary contains the points $p$ and $q$ in the respective cases (i) and (ii).   Since $\Omega^{\lambda_*}\subset\Omega$, we notice that
$$
\Sigma \subset \left(\mathbb R^N\setminus\overline{\Omega}\right)\cap\{ x \in \mathbb R^N : x\cdot\gamma < \lambda_*\}\subset D.
$$

Let $x^{\lambda_*}$ denote the reflection of a point $x\in\mathbb R^N$ with respect to $\pi_{\lambda_*}$, namely,
\begin{equation}\label{xlambda}
x^{\lambda_*}= x +2[\lambda_*-(x\cdot\gamma)]\gamma.
\end{equation}
Using the functions $v^\pm$ defined in \eqref{inside and outside}, we introduce the functions $w^\pm=w^\pm(x)$ by
\begin{equation}
\label{inside and outside reflected}
\begin{aligned}
w^+(x) &:=v^+(x)-v^+(x^{\lambda_*}) \quad \mbox{for } x \in \overline{\Omega^{\lambda_*}}, \\
w^-(x) &:= v^-(x)- v^-(x^{\lambda_*})\quad \mbox{for } x \in \overline{\Sigma}.
\end{aligned}
\end{equation}
It then follows  from \eqref{interior}--\eqref{at infinity} that
\begin{eqnarray}
& -\sigma_+\Delta w^++w^+=0 \ \mbox{ in }\Omega^{\lambda_*} &\mbox{and }\  w^+ \ge 0\ \mbox{ on }\partial \Omega^{\lambda_*},\label{eqs of w+}\\
&-\sigma_-\Delta w^-+w^-=0 \ \mbox{ in } \Sigma &\mbox{and }\  w^- \ge 0\ \mbox{ on }\partial\Sigma, \label{eqs of w-}
\end{eqnarray}
and hence by the  maximum principle
\begin{equation}
\label{positivity of w+-}
w^+ \ge 0\ \mbox{ in }\Omega^{\lambda_*} \quad \mbox{and} \quad w^- > 0 \ \mbox{ in } \Sigma.
\end{equation}
Note that $w^+$ can be zero  in $\Omega^{\lambda_*}$ since some connected component $\Omega_j$ of $\Omega$ can be symmetric with respect to $\pi_{\lambda_*}$ and, in such a case, $w^+ \equiv 0$ in $\Omega_j$. But $w^-$ is strictly positive in $\Sigma$ since $\Omega$ is not symmetric with respect to $\pi_{\lambda_*}$.

Let us first consider  the case (i).  The first equality in \eqref{transmission} yields that $w^+(p)=w^-(p)=0$. Then, it follows from  \eqref{positivity of w+-} and Hopf's boundary point lemma  that
\begin{equation}
\label{Hopf applies to case (i)}
\frac {\partial w^+}{\partial \nu}(p) \le 0 < \frac {\partial w^-}{\partial \nu}(p),
\end{equation}
where we used the fact that $\nu$ is the  outward unit normal vector to $\partial\Omega$ as well as  the inward unit normal vector to $\partial\Sigma$. It thus follows from the definition \eqref{inside and outside reflected} of $w^\pm$ that
$$
\frac{\partial v^+(x)}{\partial \nu}\bigg\vert_{x=p} \le \frac{\partial (v^+(x^{\lambda_*}))}{\partial \nu}\bigg\vert_{x=p} \quad \mbox{and} \quad
\frac{\partial v^-(x)}{\partial \nu}\bigg\vert_{x=p} > \frac{\partial (v^-(x^{\lambda_*}))}{\partial \nu}\bigg\vert_{x=p}.
$$
Reflection symmetry with respect to the plane $\pi_{\lambda_*}$ yields that
\begin{equation}
\label{result coming from symmetry}
\frac{\partial (v^\pm (x^{\lambda_*}))}{\partial \nu}\bigg\vert_{x=p} = \frac {\partial v^\pm}{\partial \nu}(p^{\lambda_*}).
\end{equation}
Indeed,  we observe that
$$
\nu(p)\cdot\gamma = - \nu(p^{\lambda_*})\cdot\gamma\ \mbox{ and }\  \nu(p) - (\nu(p)\cdot\gamma)\gamma=\nu(p^{\lambda_*}) - (\nu(p^{\lambda_*})\cdot\gamma)\gamma,
$$
and by using \eqref{xlambda}, we see that
$$
\nabla (v^\pm(x^{\lambda_*}))=(\nabla v^\pm)(x^{\lambda_*}) - 2 \left( (\nabla v^\pm)(x^{\lambda_*}) \cdot \gamma \right) \gamma.
$$
Then, combing these equalities yields \eqref{result coming from symmetry}.
It thus follows that
\begin{equation}
\label{by the reflection}
\frac {\partial v^+}{\partial \nu}(p) \le \frac {\partial v^+}{\partial \nu}(p^{\lambda_*})\quad \mbox{and} \quad \frac {\partial v^-}{\partial \nu}(p) > \frac {\partial v^-}{\partial \nu}(p^{\lambda_*}).
\end{equation}
On the other hand, the second equality in \eqref{transmission} shows that
$$
\sigma_+\frac {\partial v^+}{\partial \nu\ }(p)=  \sigma_-\frac {\partial v^-}{\partial \nu\ }(p)\quad \mbox{and} \quad \sigma_+\frac {\partial v^+}{\partial \nu\ }(p^{\lambda_*})=  \sigma_-\frac {\partial v^-}{\partial \nu\ }(p^{\lambda_*}),
$$
which contradict \eqref{by the reflection}.

Let us proceed to the case (ii).  As in \cite{Se1971},  by a translation and a rotation of coordinates, we may assume:
$$
 \gamma=(1,0,\dots,0), \ \ q=0, \ \ \lambda_*=0 \ \mbox{ and } \ \nu(q)= (0,\dots,0,1).
$$
Since $\partial\Omega$ is of class $C^2$, there exists a $C^2$ function $\varphi : \mathbb R^{N-1}\to\mathbb R$ such that
in a neighborhood of $q=0$, $\partial\Omega$ is represented as a graph $x_N=\varphi(\hat{x})$ where $\hat{x}=(x_1,\dots,x_{N-1})\in\mathbb R^{N-1}$, where
$$
\varphi(0)=0, \quad \nabla\varphi(0)=0, \quad \mbox{ and} \quad \nu=\frac 1{\sqrt{1+|\nabla\varphi|^2}}(-\nabla\varphi,1).
$$
Since the event (ii) occurs at $\lambda =0$, we observe that the function $\frac {\partial\varphi}{\partial x_1}(0,x_2,\dots,x_{N-1})$ achieves its local maximum $0$ at $(x_2,\dots,x_{N-1})=0\in \mathbb R^{N-2}$, and hence
\begin{equation}
\label{second derivatives vanish}
\frac {\partial^2 \varphi}{\partial x_1\partial x_j}(0)=0\ \mbox{ for } j=2,\dots, N-1.
\end{equation}
Notice that
\begin{equation}
\label{w+- functions}
w^\pm(x)=v^\pm(x_1,x_2,\dots,x_N)-v^\pm(-x_1,x_2,\dots,x_N),
\end{equation}
since $x^{\lambda_*}=(-x_1,x_2,\dots,x_N)$.

The equalities \eqref{transmission} at  $(\hat{x},\varphi(\hat{x}))$ in a neighborhood of $q=0$ are read as
\begin{eqnarray}
\label{transmission near q=0}
&&v^\pm=a^*, \label{isothermic interface}\\
&&\sigma_+\left(-\sum_{k=1}^{N-1}\frac {\partial\varphi}{\partial x_k}\frac {\partial v^+}{\partial x_k}+\frac {\partial v^+}{\partial x_N}\right)=
\sigma_-\left(-\sum_{k=1}^{N-1}\frac {\partial\varphi}{\partial x_k}\frac {\partial v^-}{\partial x_k}+\frac {\partial v^-}{\partial x_N}\right).\label{Snell}
\end{eqnarray}
Differentiating \eqref{isothermic interface} in $x_i$ for $i=1,\dots,N-1$ yields that at  $(\hat{x},\varphi(\hat{x}))$
\begin{equation}
\label{first tangential derivatives}
\frac {\partial v^\pm}{\partial x_i} + \frac {\partial v^\pm}{\partial x_N}\frac {\partial \varphi}{\partial x_i}=0.
\end{equation}
Then, differentiating \eqref{first tangential derivatives} in $x_j$ for $j=1,\dots,N-1$ yields that at  $(\hat{x},\varphi(\hat{x}))$
\begin{equation}
\label{second tangential derivatives}
\frac {\partial^2 v^\pm}{\partial x_j\partial x_i} +\frac {\partial^2 v^\pm}{\partial x_N\partial x_i}\frac {\partial \varphi}{\partial x_j}+ \frac {\partial^2v^\pm}{\partial x_j\partial x_N}\frac {\partial \varphi}{\partial x_i}+ \frac {\partial^2 v^\pm}{\partial x_N^2}\frac {\partial \varphi}{\partial x_i}\frac {\partial \varphi}{\partial x_j}+\frac {\partial v^\pm}{\partial x_N}\frac {\partial^2 \varphi}{\partial x_j\partial x_i}=0.
\end{equation}
By letting $\hat{x}=0$ in these equalities, we obtain from \eqref{second derivatives vanish} that
\begin{equation}
\label{values of v+- at 0}
\frac {\partial v^\pm}{\partial x_i}(0)=\frac {\partial^2 v^\pm}{\partial x_1\partial x_j}(0)=0\ \mbox{ for } i=1,\dots, N-1 \mbox{ and } j=2,\dots, N-1.
\end{equation}
Next, differentiating \eqref{Snell}  in $x_i$ for $i=1,\dots,N-1$ and letting $\hat{x}=0$ give
\begin{equation}
\label{derivative of Snell}
\sigma_+\frac {\partial^2 v^+}{\partial x_i\partial x_N}(0)=\sigma_-\frac {\partial^2 v^-}{\partial x_i\partial x_N}(0)\ \mbox{ for } i=1,\dots, N-1.
\end{equation}

Since the functions $w^\pm$ are expressed as \eqref{w+- functions}, with the aid of  \eqref{values of v+-  at 0} we have that
\begin{equation}
\label{values of w+- at 0}
w^\pm(0)= \frac {\partial w^{\pm}}{\partial x_j}(0) = \frac {\partial^2 w^\pm}{\partial x_1\partial x_j}(0)=0\ \mbox{ for }   j=1,\dots, N-1.
\end{equation}
The relations \eqref{eqs of w+}--\eqref{positivity of w+-} enable us to apply Serrin's corner point lemma (see \cite[Lemma S, p. 214]{GNN1979} or \cite[Serrin's Corner Lemma, p. 393]{R1997}) to show that
\begin{equation}
\label{second derivatives}
\frac {\partial^2 w^+}{\partial s_{+}^2}(0) \ge 0\ \mbox{ and }\ \frac {\partial^2 w^-}{\partial s_{-}^2}(0) > 0\ \mbox{ with } s_{\pm}=-\gamma\mp\nu =(-1,0,\dots,0,\mp1),
\end{equation}
where $\frac {\partial^2 w^\pm}{\partial s_{\pm}^2}$ denotes the second derivative of $w^\pm$ in the direction of $s_{\pm}$. Note that each of the directions $s_\pm$ respectively enters $\Omega^{\lambda_*}, \Sigma$, transversally to both of the hypersurfaces $\partial\Omega$ and $\pi_{\lambda_*}$.
Thus, we have from \eqref{w+- functions} and \eqref{values of w+- at 0} that
\begin{equation}
\label{strictly positive second derivatives}
 \frac {\partial^2 w^\pm}{\partial s_{\pm}^2}(0)=\pm2 \frac {\partial^2 w^\pm}{\partial x_1\partial x_N}(0)=\pm4\frac {\partial^2 v^\pm}{\partial x_1\partial x_N}(0).
\end{equation}
It then follows from \eqref{second derivatives} that
\begin{equation}
\label{second Hopf}
\frac {\partial^2 v^-}{\partial x_1\partial x_N}(0) < 0 \le \frac {\partial^2 v^+}{\partial x_1\partial x_N}(0),
\end{equation}
which contradicts \eqref{derivative of Snell} with $i=1$.
Thus $\Omega$ is symmetric with respect to $\pi_{\lambda_*}$. Since the unit vector $\gamma$ is arbitrary, $\Omega$ must be a ball and Theorem \ref{symmetry theorem} is proved.


\end{document}